\documentclass[12pt]{article}
\usepackage{times}
\usepackage{mathptm}
\usepackage{amsfonts}
\usepackage{amsmath}
\usepackage{dsfont}
\usepackage{amssymb}
\usepackage{theorem}
\usepackage{pifont}
\usepackage{mathrsfs}
\usepackage[dvips]{color}
\usepackage{xcolor}
\usepackage{graphics}
\usepackage{epsfig}

\usepackage{leqno}
\newcommand{\D}{\Delta}

\newcommand{\no}{\nonumber}
\newcommand{\n}{\nabla}

\newcommand{\p}{\partial}
\newcommand{\eps}{\varepsilon}

\newcommand{\cF}{{\mathcal{F}}}

\newcommand{\R}{\mathbb{R}}

\title{The statistical theory of the angiogenesis equations}
\author{Bj\"orn Birnir,\\
CNLS
and Department of Mathematics\\
UC Santa Barbara,\\ 
Luis Bonilla, Manuel Carretero and Filippo Terragni\\
G. Mill\'an Institute and Department of Mathematics\\
Universidad de Carlos III, Madrid
}
\date{\today}
\begin{document}
\maketitle

\begin{abstract}
Angiogenesis is a multiscale process by which a primary blood vessel issues secondary vessel sprouts that reach regions lacking oxygen. Angiogenesis can be a natural process of organ growth and development or a pathological induced by a cancerous tumor.  A mean field approximation for a stochastic model of angiogenesis consists of partial differential equation (PDE) for the density of active tip vessels. Addition of Gaussian and jump noise terms to this equation produces a stochastic PDE that  defines an infinite dimensional L\'evy process and is the basis of a statistical theory of angiogenesis. The associated functional equation has been solved and the invariant measure obtained. The results are compared to a direct numerical simulation of the stochastic model of angiogenesis and invariant measure multiplied by an exponentially decaying factor.  The results of this theory are compared to direct numerical simulations of the underlying angiogenesis model. The invariant measure and the moments are functions of the Korteweg-de Vries soliton which approximates the deterministic density of active vessel tips.
\end{abstract}

\section{Introduction}
Angiogenesis is the process of cells organizing themselves into into blood vessels that grow from existing vessels and 
carry blood to organs and through tissue. It occurs in normal conditions of organ growth and regeneration, wound healing and tissue repair, and also in pathological conditions such as cancer, diabetes, rheumatoid arthritis or neovascular age-related macular degeneration. 

Angiogenesis is driven by Vessel Endothelial Growth Factor (VEGF) and other pro-angiogenic proteins which are secreted by cells experiencing lack of oxygen (hypoxia). VEGF diffuses in the tissue, binds to extracellular matrix (ECM) components and forms a spatial concentration gradient in the direction of hypoxia. Once VEGF molecules reach an existing blood vessel, the latter walls open as a response and new vessel sprouts grow out of endothelial cells (ECs) off the vessel. Through the cellular signaling Notch process, VEGF activates the tip cell phenotype in ECs, which then grow filopodia with many VEGF receptors. The tip cells pull the other ECs (called stalk cells), open a pathway in the ECM, lead the new sprouts, and migrate in the direction of increasing VEGF concentration \cite{ger03}. Signaling and mechanical cues between neighboring ECs cause {\em branching of new sprouts} \cite{hel07,jol15,veg20}. Stalk cells in growing sprouts alter their shape to form a lumen (wall of the sprout) connected to the initial vessel that is capable of carrying blood \cite{geb16}. Sprouts meet and merge in a process called {\em anastomosis} to improve blood circulation inside the new vessels. Poorly perfused vessels may become thinner and their ECs, in a process that inverts angiogenesis, may retract to neighboring vessels leading to a more robust blood circulation \cite{fra15}. Thus, the vascular plexus remodels into a highly organized and hierarchical network of larger vessels ramifying into smaller ones \cite{szy18}. In normal processes of wound healing or organ growth, the cells inhibit the production of growth factors when the process is finished. In pathological angiogenesis, e.g., cancer, tumor cells lack oxygen and nutrients and produce VEGF that induces angiogenesis from a nearby primary blood vessel. The generated new vessel sprouts move and reach the tumor \cite{fol74,car05}. Tumor cells continue secreting growth factors that attract more vessel sprouts and facilitate their expansion. 

Together with experiments, many models spanning from the cellular to macroscopic scales try to understand angiogenesis; see the reviews \cite{and98,bon19,byr10,hec15,man04,pla04,sci13,vil17}. Early models consider reaction-diffusion equations for densities of cells and chemicals (growth factors, fibronectin, etc.) \cite{lio77,byr95,cha93,cha95}. They cannot treat the growth and evolution of individual blood vessels. Tip cell stochastic models of tumor induced angiogenesis are among the simplest ones for this complex multiscale process. Their basic assumptions are that (i) the cells of a blood sprout tip do not proliferate and move towards the tumor producing growth factor, and (ii) proliferating stalk cells build the sprout along the trajectory of the sprout tip. Thus tip cell models are based on the motion of single particles representing the tip cells and their trajectories constitute the advancing blood vessels network \cite{sto91,sto91b,and98,pla04,man04,BCAC14,hec15,bon19,ter16}. Tip cell models describe angiogenesis over distances that are large compared with a cell size, thereby not incorporating descriptions of cellular and subcellular scales. These models are typically random: the motion of each tip typically includes directed Brownian motion, branching and anastomosis are birth and death processes, respectively. The random evolution of the tip cells may affect and be affected by reaction-diffusion equations for VEGF and other quantities (hybrid models). Other models contemplate motion of the tip cells on a lattice \cite{and98,pla04}, and are related to cellular automata models \cite{pil17,mar20,mar21}. More complex models include tip and stalk cell dynamics, the motion of tip and stalk cells on the extracellular matrix outside blood vessels, shape and size changes of cells, signaling pathways and EC phenotype selection, blood circulation in newly formed vessels, and so on \cite{bau07,jac10,tra11,sci13,ben14,oer14,hec15,jol15,per17,vil17,ber18,bon19,veg20}. The evolution of the blood vessels governed by these processes may be difficult to measure but various methods including laboratory experiments can be used to measure the statistical properties of an angiogenic network. These measurements of statistical quantities of the network can then be compared to simulations.

The authors have analyzed a system of stochastic differential equations plus birth and death processes, from \cite{BCAC14,ter16,bon20}, describing vessel growth and the associated PDE describing the evolution of the density of active vessel tips. In \cite{BB16}, they found that the growth of the tip of the blood vessel is well described by the KdV soliton. This leads to a control theory of angiogenesis that is currently being developed. In this paper, we will add the noise associated with branching and anastomosis, to the density equations, and develop the statistical theory of angiogenesis. Not surprisingly the solitons are found to play a major role in this theory and the noise in both branching and anastomosis is found to be multiplicative and depend only on the density.

The paper is organized in the following matter. In Section 2, we derive the density equations and formulated the noise in them. In Section 3, we explain the log-Poisson process arising from the noise in anastomosis and derive the geometric Brownian motion resulting from the noise in branching. In Section 4, we explain the statistical theory by the example of the stochastic Navier-Stokes Equation and derive the invariant measure of the stochastic angiogenesis equation. In Section 5, we compute the moments of the density and compare with simulations. In Section 6, we compute the structure function of the density. Section 7, contains a discussion. In Appendix A, we include a short summary of jumps and L\'evy processes In Appendix B, we indicate how to calculate the moments of the density from numerical simulations of the angiogenesis stochastic process.

\section{The Density Equation}
We start with the equation for the density of active tip cells $p$ in nondimensional form
\begin{equation}
\label{eq:angio}
\frac{\p p}{\p t} =\!\left[ \frac{AC}{1+C} \ \delta_v(v-v_0)- \Gamma S\right]\!p- v\cdot \n_x p-\n_v\left[
\left(\frac{\delta_1  \n_x C}{(1+\Gamma_1C)^{q}}-\sigma v\right)\!p\right]+\frac{\sigma}{2}\D_v p
\end{equation}
from \cite{BCAC14,ter16,bon20}, where
\[
\tilde p(t,x) = \int p(t,x,v)dv,\quad S(t,x)= \int_0^t \tilde p(s,x)ds,
\]
are the marginal density of $p(t,x,v)$ and the density of stalk cells, respectively. $\delta_v(v)$ is a Gaussian function having zero mean and small variance $\sigma_v$. The density equation is derived from the Langevin equations for the blood vessel extension, see \cite{BCAC14,ter16,bon20}:
\begin{eqnarray}
dX^k(t)&=&v^k(t)dt,\nonumber\\
dv^k(t)&=&[-\sigma v^k(t)+F(C(t,X^k(t)))]dt+\sqrt{\sigma} dW^k(t),\label{eq2}
\end{eqnarray}
for $t>T^k$, the random time when the $k$th tip, located at $X^k(t)$ and moving with velocity $v^k(t)$ appears as a consequence of another tip's branching. When an active tip arrives at a point that was occupied by another tip at a previous time, it disappears, which corresponds to {\em anastomosis} or loop formation \cite{BCAC14,ter16,bon20}. $W^k(t)$ denotes independent Brownian motion and the derivation uses Ito's lemma, where the Brownian noise $\sqrt{\sigma} dW^k(t)$ produces the Laplacian term $\frac{\sigma}{2}\D_v p$ in the equation. The chemotactic force is
\begin{eqnarray}
F(C) = \frac{ \delta_1}{(1+\Gamma_1C)^q} \nabla_xC,\label{eq3}
\end{eqnarray}
where $C$ is the VEGF concentration (called tumor angiogenic factor or TAF in tumor induced angiogenesis). In the hybrid stochastic model, the equation for the TAF density $C(t,x)$ involves diffusion and consumption by the advancing tip cells \cite{bon20},
\begin{eqnarray} 
\frac{\partial}{\partial t}C(t,x)\!&\!=\!& \! \kappa \Delta_x C(t,x)
-\chi C(t,x) \sum_{i=1}^{N(t)} |v^i(t)|\, \delta_{\sigma_x}(x-X^i(t)),  \label{eq4}
\end{eqnarray}
where $N(t)$ is the number of active tips at time $t$ and $\delta_{\sigma_x}$ are regularized delta functions:
\begin{eqnarray}
\delta_{\sigma_x}(\mathbf{x})= \frac{e^{-x^2/\sigma_x^2}\, e^{-y^2/\sigma_y^2}}{\pi\sigma_x \sigma_y}, \quad \mathbf{x}=(x,y). \label{eq5}
\end{eqnarray}
The source terms on the right hand side (RHS) of Eq.~\eqref{eq:angio} arise from branching and anastomosis of active tips. The probability that a tip branches from one of the existing ones during an infinitesimal time interval $(t, t + dt]$ is proportional to $\sum_{i=1}^{N(t)}\alpha(C(t,X^i(t)))dt$, where  
\begin{eqnarray}
\alpha(C)=\frac{A\, C}{1+C},\quad A>0. \label{eq6}
\end{eqnarray}
The velocity of the new tip that branches from tip $i$ at time $T^i$ is selected out of a normal distribution, $\delta_{\sigma_v}(v-v_0)$, with mean $v_0$ and a narrow variance $\sigma_v^2$. The regularized delta function $\delta_{\sigma_v}(\mathbf{x})$ is given by Eq.~\eqref{eq5} with $\sigma_x=\sigma_y=\sigma_v$. 

When using the deterministic description of Eq.~\eqref{eq:angio}, the $C$ the tumor angiogenic factor satisfies the diffusive mean-field equation
\begin{eqnarray}
\frac{\p }{\p t} C(t,x) = \kappa \Delta_x C(t,x)-\chi C(t,x)\, j(t,x);\label{eq7}
\end{eqnarray}
instead of Eq.~\eqref{eq4}. Here 
\[
j(t,x) = \int |v|\, p(t,x,v)dv
\]
is the current density, see \cite{bon20}. Representative values of all involved dimensionless parameters can be found in Refs.~\cite{BCAC14,ter16,bon20}.

Density, marginal density and current density are the expected values of 
\begin{eqnarray*}
&&\sum_{i=1}^{N(t)}\delta_{\sigma_x}(x-X^i(t))\, \delta_{\sigma_v}(v-v^i(t)),\\
&&\sum_{i=1}^{N(t)}\delta_{\sigma_x}(x-X^i(t)), \\
&&\sum_{i=1}^{N(t)}\!\!|v^i(t)|\, \delta_{\sigma_x}(x-X^i(t)),
\end{eqnarray*}
respectively, with respect to the stochastic process, as the regularizations of the delta functions disappear \cite{ter16,bon20}. Equations \eqref{eq:angio} and \eqref{eq7} are mean field approximations for these expected values and for the average TAF concentration.  Moments of the density can be directly calculated from simulations of the stochastic process, as we shall see later. However, we want to build a theory of these moments by adding appropriate noise terms to the equation \eqref{eq:angio} and then analyzing the resulting stochastic PDE. The first term on the RHS of the equation (\ref{eq:angio}) represents a (multiplicative) jump term and we add a multiplicative jump noise term
\[
p\int_{\mathbb{R}}h(z,v,t){\bar N}(dz,dt)
\]
associated with such jumps. The justification for this noise term is that tip branching and anastomosis are not deterministic processes, although they are represented in the density equation above as deterministic terms, in the mean field approximations. Thus we can expect noise to be associated to these terms and since both trip branching and anastomosis create jumps in the density, the correct form of the associated noise is multiplicative Poissonian jump noise. This gives the equation
\begin{eqnarray}
\label{eq:angiojumps}
\frac{\p p}{\p t} =[ \alpha(C) \ \delta_v(v-v_0)- \Gamma  S]\, p- v\cdot \n_x p-\n_v\cdot\!\left[
(F(C)-\sigma v)p\right]\\
+ \eps\left[\sum_{k\neq 0}\left(d_k+c_k^{1/2}\frac{db_t^k}{dt}\right) <e^{ik\cdot v},p>\right] p+\frac{\sigma}{2}\D_v p+\eps p\int_{\mathbb{R}}h(z,v,t){\bar N}(dz,1),\no
\end{eqnarray}
where $\alpha(C)=\frac{AC}{1+C}$ and $b_t^k$ are i.i.d. Brownian motions. The first noise term is associated to velocity fluctuations through the Fourier coefficients $<e^{ik\cdot v},p>$. Here $h(z,v,t)=h_b(z,v,t)+h_a(z,t)$. We expect the noise to be small, $\eps << 1$, in many cases but we will allow it to be as large as $\eps =1$, in the computations below.

\[
h_a(z,t)=-\Gamma S,
\] 
and
\[
h_b(z,v,t)= \alpha(C) \delta_v(v-v_0),
\] 
denote the sizes of the jumps, associated with anastomosis and branching respectively, and ${\bar N}$ is the compensated number (of jumps) process, see \cite{OS05}. 

The Fourier coefficients $d_k$ and $c_k^{1/2}$ are summable, $\sum_{k\neq 0} |d_k| < \infty$, $\sum_{k\neq 0} |c_k^{1/2}| < \infty$. The $c_k^{1/2}$ are also square summable, $\sum_{k\neq 0} c_k < \infty$. Thus the first and the last terms in the second line in the equation represent respectively the continuous (Brownian) and discrete (jump) noise. The Brownian terms are accompanied by  a deterministic estimate for the large deviation (the $d_k$s). We have made the Brownian terms as generic as possible by modeling noise in every (velocity) Fourier component of $p$. The stochastic partial differential equation (\ref{eq:angiojumps}) must be accompanied by vanishing boundary conditions on the plane  $x \in \mathbb{R}^2$, and periodic boundary conditions in $v \in \mathbb{T}^3$. An initial condition for the density must also be specified. For the existence theory in the deterministic case, see \cite{CD16} and \cite{CDN17}, for convergence of positivity preserving numerical schemes, see \cite{BCCDNT18}.

\section{The Log-Poisson Processes}

An integrating factor can be found to simplify the density equation (\ref{eq:angiojumps}). It uses the Ito-L\'evy formula and the geometric L\'evy process. Let $f(X^k,v^k) = \ln(q_P)$ and consider Ito-L\'evy's formula (\ref{eq:Ito-Levy formula}) in Appendix A, 
\[
d\ln(q_P) =\int_{\mathbb{R}}[\ln(1+h(t,z))-h(t,z)]m(dz)dt+\int_{\mathbb{R}}[\ln(1+h(t,z))]{\bar N}(dz,dt),
\]
where we have only taken the Poisson noise into account, 
\[
dq_P = q_P\int_{\mathbb{R}}h(t,z){\bar N}(dz,dt).
\]
Then if $h=\beta-1$ and $E(N_t)=-\frac{\gamma \ln(b)}{\beta-1}$, this means that neither $h$ nor $N_t$ depend on $z$ and the integral above reduces to a product. $N_t$ is the Poisson number process and $-\frac{\gamma \ln(b)}{\beta-1}$ its mean. The two parameters $\gamma$ and $b$ will be assigned values shortly. 
\[
d\ln(q_P)= N_t \ln(\beta) - (\beta-1)E(N_t)= N_t \ln(\beta)+(\beta-1)\frac{\gamma \ln(b)}{\beta-1}=N_t \ln(\beta)+\gamma \ln(b),
\]
see \cite{BB13}.  Provided $q_P(0)=1$, solving for $q$ gives
\[
q_P = b^\gamma \beta^{N_t}.
\]
This log-Poisson process is the integrating factor that we will use below. In \cite{BB13}, Example 1.5, it is shown how to compute the moments of a log-Poisson processes: %
\begin{enumerate}
\item
Suppose the Poisson process $N_k$ has the mean $\lambda = -\frac{\gamma \ln|k|}{\beta-1}$ (i.e., $b=|k|$ above), $k$ is the wave number, then it is straight-forward to compute the mean of the log-Poisson process $|k|^\gamma \beta^{N_k}$. Namely, 
\[
E(|k|^\gamma \beta^{N_k})=\sum_{j=0}^\infty |k|^\gamma \beta^j \frac{\lambda^j}{j!} e^{-\lambda}=|k|^\gamma \sum_{j=0}^\infty \frac{(\beta\lambda)^j}{j!} e^{-\lambda}=|k|^\gamma e^{(\beta-1)\lambda}.
\]
Thus
\[
\ln[E(|k|^\gamma \beta^{N_k})]=
\gamma \ln|k|+(\beta-1)\lambda=\gamma \ln|k|-(\beta-1)\frac{\gamma}{\beta-1}\ln|k| =0,
\]
and we get the mean 
\begin{eqnarray}
\label{eq:mean}
E(|k|^\gamma \beta^{N_k})=1.
\end{eqnarray}
\item 
Now we compute the $n$th moment  $E([|k|^\gamma \beta^{N_k}]^{n})$ of the log-Poisson process $|k|^\gamma \beta^{N_k}$ above.
By a similar computations as above, 
\[
E([|k|^\gamma \beta^{N_k}]^{n})=|k|^{n\gamma} e^{(\beta^{n}-1)\lambda},
\]
where $\lambda$ is the mean from Part 1, 
therefore,
\[
\ln[E([|k|^\gamma \beta^{N_k}]^{n})]=n\gamma \ln|k|+(\beta^{n}-1)\lambda=(n-\frac{(\beta^{n}-1)}{\beta -1})\gamma \ln|k|.
\]
Finally, we get that
\begin{eqnarray}
\label{eq:nmom}
E([|k|^\gamma \beta^{N_k}]^{n})=|k|^{\gamma (n-\frac{(\beta^{n}-1)}{\beta -1})}.
\end{eqnarray}
This will be the case for the example of the stochastic Navier-Stokes Equation below, we get a log-Poisson process for each wave number. 
\item The Geometric Brownian motion gives the factor
\[
q_B=\exp\{ \sum_{k\neq 0} (d_k-\frac{1}{2}c_k)t+c_k^{1/2}b_t^k\}
\]
see formula (\ref{eq:z}) for $z$ in Appendix A, following \O ksendal \cite{Ok03},
\begin{eqnarray}
\label{eq:geob}
E(q_B)=\exp\{ \sum_{k\neq 0} d_k\ t\ + \sum_{k\neq 0}c_k^{1/2}b_0^k \}=\exp{(-d\ t-d_0)},
\end{eqnarray}
where the initial conditions for the Brownian motions,  starting at zero, are $b_0^k$. (Ito's formula produces the term $(1/2)c_k$ from the Brownian motion $b_t^k$ that cancels the term $-(1/2)c_k$ in the corresponding large deviation.) $-d = \sum_{k\neq 0} d_k$ is the drift coefficient, that we take to be negative. The initial condition $d_0=-\sum_{k\neq 0}c_k^{1/2}b_0^k$, can either be positive or negative. 
\end{enumerate}

\section{The Invariant Measure}

The PDE (\ref{eq:angiojumps}) can be used to define an infinite dimensional L\'evy process. The statistical properties such as the mean and moments of this process are determined by the invariant measure on function space associated with this process. We now review the theory of the Kolmogorov-Hopf equations determining the invariant measure. The first such equation that is a functional differential equation, where the derivatives are with respect to a function in a Banach space, was written down by Hopf \cite{Ho52},
\[
\frac{\p \phi(u)}{\p t}= \langle Au, \nabla_u \phi(u) \rangle,
\]
here $\phi$ is a bounded function of the L\'evy process $u$, and $Au$ is the deterministic part of the evolution equation determining $u$, $\langle \cdot, \cdot \rangle$ is the dual pairing in the Banach space. Hopf actually worked with the equation for the characteristic function, that is equivalent to the above equation, and he was looking for the invariant measure of the Navier-Stokes equation. Hopf's equation was reportedly solved by Kolmogorov that found that the invariant measure for the deterministic Navier-Stokes equation is disappointingly $ \mu = \delta (u)$, a delta function concentrated at the origin. The reason for this is that the Navier-Stokes equation is dissipative and all solutions eventually decay to the origin. 

Da Prato and Zabczyk developed \cite{DPZ96} a method that can be used to find the invariant measure for stochastic partial differential equation (SPDE) of the form
\[
du = Au\, dt + \sum_{j \in \mathbb{Z}^n}c_j^{1/2}db^j_t e^{2\pi j \cdot x},
\]
where $A$ is an operator on the Hilbert space on an $n$-torus $L^2(\mathbb{T}^n)$.
Here the $c_j^{1/2}$ are $n$-vectors, their inner products converge $\sum_{j \in \mathbb{Z}^n} c_j \leq \infty$, $c_j =c_j^{1/2}\cdot c_j^{1/2}$, and the $b^j_t$ are independent Brownian motion accompanying each Fourier coefficient. The corresponding Kolmogorov-Hopf equation is 
\[
\frac{\p \phi(u)}{\p t}= \frac{1}{2} 
\mathrm{Tr}[C\Delta_u\phi(u)]+\langle Au, \nabla_u \phi(u) \rangle,
\]
where $C$ is the trace class matrix having the $c_j$s along the diagonal. The invariant measure producing the solution of this equation is Gaussian,
\[
\mu = \mathcal{N}_{(0,Q)}(u),
\]
where $Q= \int_0^\infty e^{tA}Ce^{tA^*}dt$ is the variance, see \cite{DP06}. 

In \cite{BB211} the invariant measure of the Navier-Stokes equation was computed using 
Stochastic Closure Theory, were the Navier-Stokes equation is split into a Reynolds Averages Navier-Stokes (RANS) equation and a stochastic Navier-Stokes (SNS) equation for the small scale flow. This is the equation 
\begin{eqnarray}
\label{eq:SNS}
du=\!\left[\nu\Delta u-u\cdot \nabla u+\nabla \Delta^{-1}\mathrm{Tr}(\nabla u)^2\right]dt  + \sum_{j \in \mathbb{Z}^n}(|k|^{1/3}d_kdt+c_j^{1/2}db^j_t) e^{2\pi j \cdot x}\no\\+\sum_{k=0}^m\int_{\mathbb{R}}h_k{\bar N}^k(dt,dz)u,
\end{eqnarray}
where we have eliminated the pressure, using the incompressibility condition and the $d_k$s terms are contributed by the large deviation principle, see \cite{BB211}. This equation can be written in the form 
\[
du = Au\, dt+\sum_{j \in \mathbb{Z}^n}(|k|^{1/3}d_kdt+c_j^{1/2}db^j_t) e^{2\pi j \cdot x}+\sum_{k=0}^m\int_{\mathbb{R}}h_k{\bar N}^k(dt,dz)u,
\]
The large deviation term in the noise will simply contribute a non-trivial mean to 
the infinite-dimensional Gaussian, but the last term that is multiplicative in $u$ has to be treated differently. 

In \cite{BB211} the invariant measure of the stochastic Navier-Stokes equation was computed by applying Girsanov's theorem. This amounts to choosing the parameters $\gamma = 2/3 = \beta$ and introducing an integrating 
factor 
\[
P_t =\prod_{k=0}^m (|k|^{2/3}(2/3)^{N^k_t})^{1/3}
\]
based on the log-Poisson processes in the last section. This means that the Kolmogorov-Hopf equation becomes
\begin{eqnarray}
\label{eq:KH1}
\frac{\p \phi(u)}{\p t}= \frac{1}{2} 
\mathrm{Tr}[P_tCP^*_t\Delta_u\phi(u)]+\langle P_t{\tilde D}, \nabla_u \phi(u) \rangle+\langle A(u)P_t, \nabla_u \phi(u) \rangle,
\end{eqnarray}
where $C$ is the trace class matrix having the $c_j$s along the diagonal, and ${\tilde D}=[|k|^{1/3}d_k]$ is an infinite vector. The invariant measure producing the solution of this equation is Gaussian,
\begin{eqnarray}
\label{eq:invarmeasure}
\mu = \mathcal{N}_{(E,Q)}\ast\mathbb{P}_{N_\infty}(u),
\end{eqnarray}
where $Q= \int_0^\infty e^{tA}P_tCP^*_te^{tA^*}dt$ is the variance, and the mean is 
$\mathrm{Tr}(EI)$, where $E=\int_0^\infty e^{tA}P_t{\bar D}dt$. ${\bar D}$ is the matrix with $|k|^{1/3}d_k$s along the diagonal and 
$I$ is the matrix with Fourier components $e_k(x)$ on the diagonal. 
$N^k_t$ in the integrating factor $P_t$ above, is the log-Poisson process and its law is:
\begin{eqnarray}
\label{eq:measure}
\mathbb{P}_{N^k_t}(\cdot)=\sum_{j=0}^\infty (\cdot)\frac{(m^k_t)^j e^{-m^k_t}}{j!} \delta_{(N^k_t-j)},
\end{eqnarray}
where $m^k_t = -\frac{\gamma \ln(|k|)}{\beta-1}$ is the mean of the Log-Poisson process for the $k$th wavenumber and $\mathcal{P}_j = \frac{(m^k_t)^j e^{-m^k_t}}{j!} $ is the probability of having exactly $j$ jumps in the $k$th wavenumber.
Thus we read the normal law with variance $Q$ from the first and the last term of the 
Kolmogorov-Hopf equation, the mean from the second term and the Poisson law from the 
integrating factor. Then one can prove directly that (\ref{eq:invarmeasure}) is the invariant measure of (\ref{eq:SNS}) and provides the solution of the Kolmogorov-Hopf equation (\ref{eq:KH1}), see \cite{BB13}. 

We will now use the same reasoning to find the Kolmogorov-Hopf equation associated to 
(\ref{eq:angiojumps}) and solve it. First we write (\ref{eq:angiojumps}) as a stochastic partial differential equation,
\begin{eqnarray}
\label{eq:sangiojumps}
dp =\left(\frac{\sigma}{2}\D_v p+ [\alpha(C) \ \delta_v(v-v_0)- \Gamma  S]\, p- v\cdot \n_x p-\n_v\cdot\!\left[
(F(C)-\sigma v)p\right]\right)dt\no\\
+\eps \left[\sum_{k\neq 0}\left(d_kdt+c_k^{1/2}db_t^k\right) {\hat p_k}\right] p+\eps p\int_{\mathbb{R}}h(z,v,t){\bar N}(dz,dt),
\end{eqnarray}
where $\hat p_k =<p,e^{ik\cdot v}>$, models velocity fluctuations. This shows that (\ref{eq:sangiojumps}) is different from the stochastic Navier-Stokes equation above in that it contains two multiplicative noise terms and no additive noise. We write it in the form
\[
dp = A(p)dt + \left(\left[\sum_{k\neq 0}\left(d_kdt+c_k^{1/2}db_t^k\right) {\hat p_k}\right] +\int_{\mathbb{R}}h(z,v,t){\bar N}(dz,dt)\right)p,
\]
where we have set $\eps =1$. Now we find an integrating factor that is a product of two terms,
\[
q_P= |p_t|^{\gamma}\beta^{N_t},
\]
a log-Poisson process as in the previous section, with the mean $\lambda=m_t = -\frac{\gamma \ln(p_t)}{\beta-1}$, and the geometric Brownian motion
\[
q_B=\exp\{ \sum_{k\neq 0} (d_k-\frac{1}{2}c_k){\hat p_k}t+c_k^{1/2}{\hat p_k}b_t^k\}
\]
as in Equation (\ref{eq:geob}). By Ito's formula
\[
E(q_B)=\exp\{ \sum_{k\neq 0} {\hat p_k}d_k\ t\ + \sum_{k\neq 0}{\hat p_k}c_k^{1/2}b_0^k \}=\exp{(-d\ t-d_0)}.
\]
Notice that the log-Poissonian is different from the stochastic Navier-Stokes equation where the jumps depend on the Fourier coefficients ($k$) of the solution, whereas for (\ref{eq:sangiojumps}) the jumps depends on the whole solution $p_t$, or are the same for all Fourier coefficients. 
Then the integrating factor becomes
\[
P_t =q_Bq_P
\]
and the Kolmogorov-Hopf equation reduces to Hopf's equation with an integrating factor
\begin{eqnarray}
\label{eq:KH2}
\frac{\p \phi(p)}{\p t}= \langle A(p)P_t, \nabla_p \phi(u) \rangle.
\end{eqnarray}
The invariant measure of (\ref{eq:sangiojumps}) is now a convolution of the Poisson distribution
with a Gaussian,
\[
\mu=\mathcal{N}_{(e,q)}\ast\mathbb{P}_{N_\infty}(p)
\]
where the mean of the Gaussian is $e_t=\sum_{k\neq 0} (d_k-\frac{1}{2}c_k){\hat p_k}$ and the variance is 
$q_t =\sum_{k\neq 0}c_k{\hat p}_k^2$. The Gaussian is the distribution of the sum of the Fourier coefficients of $p$ with respect to the velocity $v$, with mean $e_t$ and variance $q_t$. 

The solution of the equation (\ref{eq:KH2}) can be written, see \cite{BB13}, as
\begin{eqnarray}
\label{eq:KHsol}
R_t \phi&=& \frac{1}{\sqrt{\pi q_t}}\int_{-\infty}^\infty\int_0^\infty  \phi(p)\ \mathbb{P}_{N_t}e^{-\frac{(p-e_t)^2}{q_t}}dxdp\no\\
&=&\frac{1}{\sqrt{\pi q_t}}\int_{-\infty}^\infty\sum_{j=0}^\infty \int_0^\infty   \phi(p) \frac{m_t^j e^{-m_t}}{j!}\delta(x-j)e^{-\frac{(p-e_t)^2}{q_t}}dxdp\no\\&=&\frac{1}{\sqrt{\pi q_t}}\int_{-\infty}^\infty\sum_{j=0}^\infty   \phi(p) \frac{m_t^j e^{-m_t}}{j!}e^{-\frac{(p-e_t)^2}{q_t}}dp,
\end{eqnarray}
where $m_t = -\frac{\gamma \ln(p_t)}{\beta-1}$ is the mean of the Log-Poisson process and $\mathcal{P}_j = \frac{m_t^j e^{-m_t}}{j!} $ is the probability of having exactly $j$ jumps $N^j_\infty = N^j$.  For $\phi = p^n$, we get
\begin{eqnarray}
\label{eq:moments}
\frac{1}{\sqrt{\pi q_t}}\int_{-\infty}^\infty\sum_{j=0}^\infty \phi(e^{Kt}p_0P_t)\frac{m_t^j e^{-m_t}}{j!} \delta_{(N_t-j)}e^{-\frac{(p-e_t)^2}{q_t}}dp={p}^n_t E(q_P^nq_B^n)\no\\=e^{-ndt-nd_0}{ p}^{n(1-\gamma)+\gamma(\frac{\beta^n-1}{\beta-1})}_t,
\end{eqnarray}
by the computation of $E(q_B)$ and $E(q_P)$ in Equations (\ref{eq:geob}) and (\ref{eq:nmom}). 
We would like to take the limit $t \to \infty $, to obtain the invariant measure
\begin{eqnarray}
\label{eq:measure}
\mathbb{P}_{N_\infty}(\cdot) \ast \mathcal{N}_{(e,q)}(\cdot)\ =
\int_{-\infty}^\infty \sum_{j=0}^\infty (\cdot)\frac{m_\infty^j e^{-m_\infty}}{j!} 
\delta_{(N_\infty-j)}\frac{1}{\sqrt{\pi q_\infty}}e^{-\frac{(p-e_\infty)^2}{q_\infty}}dp,
\end{eqnarray}
where $m_\infty = -\frac{\gamma \ln(p_\infty)}{\beta-1}$ and $e= e_\infty,\: q=q_\infty$ are the limits as $t \to \infty$. However, this would give us the trival limit zero, because of the decay $e^{-ndt}$ due to the multiplicative Brownian noise.
Thus we take the limit $p_t \to p_\infty$ and consider the long time statistics to be determined by 
\begin{eqnarray}
\label{eq:expinv}
e^{-ndt-nd_0}{ p}^{n(1-\gamma)+\gamma(\frac{\beta^n-1}{\beta-1})}_\infty.
\end{eqnarray}
This can be thought of as an invariant measure multiplied by an exponentially decaying factor. It gives an excellent approximation to the real statistics for large times.

\section{Comparison with Simulations and Theory} 
In \cite{BB16}, the authors showed that the density $p$ evolves towards a 1D (Korteweg-de Vries) soliton profile for a slab geometry in which the primary vessel is the $y$ axis and the tumor is centered at $x=L$, $y=0$. A general 2D geometry is discussed in \cite{bon20}. For the 1D slab geometry, $p$ is a product of a Maxwellian distribution in velocity, centered at $v_0$, a transversal Gaussian function that approaches $\delta(y)$, and the soliton from \cite{BB16},
\begin{eqnarray}
&&p(t,x,v)=\frac{e^{-|v-v_0|^2}}{\pi}\delta(y)\,\tilde{p}(x,t), \label{eq22}\\
&&\tilde{p}(x,t)=\frac{(2K\Gamma+\mu^2)c}{2\Gamma(c-F_x)}\mbox{sech}^2\left(\frac{\sqrt{2K\Gamma+\mu^2}}{2\Gamma(c-F_x)}
(x-ct+x_0)\right)\!, \label{eq23}
\end{eqnarray}
see \cite{bon20}. The justification of the Maxwellian distribution is that the source term in Equation (\ref{eq:sangiojumps}) selects velocities in a small neighborhood of $v_0$ because they are the only velocities for which the birth term in the equation can compensate the anastomosis death term. Then a Chapman-Enskog expansion reduces Equation (\ref{eq:sangiojumps}) to the equation in \cite{BB16} that has the soliton solution, see \cite{BCT16}. The factor $\delta(y)$ follows from a multiple scales method explained in \cite{bon20}. 
It is reasonable to expect the averaged density and its moments (\ref{eq:moments}) to be functions of this soliton. This is in fact the case, but we now spell out the assumption that we make to establish this relation.
The jump associated with Poisson process is:
\[
h= h_b+h_a=\left(\mathds{1}_b\alpha(C)\,\delta_v(v-v_0)-\mathds{1}_a\frac{1}{2}\Gamma 
\frac{e^{-|v-v_0|^2}}{\pi}\int_0^t {\tilde p} ds\right)\!,
\]
where 
$\mathds{1}_a$ and $\mathds{1}_b$ are one during branching and anastomosis respectively and zero otherwise. Now substituting in the soliton for $p=a\, \mbox{sech}^2(b(x-ct+x_o))$, 
we get that
\[
S(t,x)=\int_0^t {\tilde p} ds=-(a/bc)\tanh(b(x-ct+x_o))+(a/bc)\tanh(b(x+x_o)).
\]
This expression ranges from $0$ to $2(a/bc)$, with $x$ and we choose the average $S(\infty)=(a/bc)$.  Another way of determining this limit, is to let $x+x_o=ct$ and take the limit as $x \to \infty$ and $t \to \infty$.
The (velocity) mean of the jump is then, the average of $h$ over the velocity,
\[
\overline{h} = \int_{-\infty}^\infty h(x,v,t)dv = \overline{h}_b+\overline{h}_a=(\mathds{1}_b\overline{\alpha} - \mathds{1}_a\frac{1}{2} \Gamma \frac{a}{bc}).
\]
This gives the average rate
\[
\overline{m}= -\frac{\gamma}{(\mathds{1}_b \overline{\alpha} -\mathds{1}_a\frac{1}{2} \Gamma S(\infty))}\log({\overline {\tilde p}}).
\]

The averaged branching rate $\overline m_b=24.70$ and the averaged anastomosis rate is $\overline m_a=-22.49$, see Table 1. However, the rates in Equation (\ref{eq:moments}) depend on $x$ (and $t$ once we switch on the time evolution), so we must average them to get the average rates. We use the ergodic hypothesis to do this. Strictly speaking we do not have a non-trivial invariant measure to apply the ergodic theory, as discussed in the last section. However, we can ignore the decaying part $q_B$ of the solution and use (\ref{eq:expinv}) for the purpose of computing the averaged branching and anastomosis rates.  We take the log of both sides of the equation from last section,
\[
e^{-m\frac{h}{\gamma}}=p,
\]
and approximate the time averages of both sides by the averages with respect to (\ref{eq:measure})
\[
-{\overline m}\frac{\overline{h}}{\gamma}=\frac{1}{\sqrt{\pi q}}\int_{-\infty}^\infty |\mathbb{P}_{N_\infty}(\log(p))|e^{\frac{-(p-e)^2}{q}}dp=\log(p^*)+E[\log(q_P)]=\log(p^*),
\]
since $E[\log(q_P)] \leq \log(E[q_P]) = 0$, by Jensen's inequality, and $E[q_P]=1$, by Equation (\ref{eq:mean}). 
The first equality above is by the ergodic hypothesis, using the invariant measure (\ref{eq:measure}), the second equality is by the mean value theorem, where $p^*$ is the value of $p$ as some fixed velocity $v^*$.  We have set $p=p_\infty q_P$ and have ignored the decaying part $q_B$. Finally, we approximate $p^*$ by the spatial mean of $\tilde p$
\[
\log(p^*) \approx \log(\int_{-\infty}^\infty {\tilde p} dx)= \log({\overline {\tilde p}}).
\]
A substitution of the soliton above into the integral gives 
\[
\overline m = -\frac{\gamma}{\overline{h}}\ln B,
\]
where $B= \frac{2(\sqrt{2K\Gamma+\mu^2})c}{\Gamma}$. This allows to compute the exponents
\[
\gamma_b = -\frac{\overline{\alpha}}{\ln B}\ {\overline m_b}=-14.58\:\:\: \mathrm{and} \:\:\: 
\gamma_a =  \frac{\Gamma S(\infty)}{2\ln B}\ {\overline m_a}=-9.87,
\]
in Equation (\ref{eq:moments}), using the values $K=266,\ c=3,\ \Gamma = 0.32,\ \mu=8.5$, used in the simulations, see \cite{BB16}, so $B=291.977$ 
and $\ln B = 5.677$. For consistency, we have set  $\overline{\alpha}=3.35$.
\begin{table}[h!]
\centering
{\small
\begin{tabular}{@{}lcccccccc@{}} 
\quad & \quad ${\overline m}_b$ & \,\,\,\,\,\, ${\overline h}_b$ & \,\,\,\,\,\,\, $\overline{m}_a$ & \,\,\,\,\,\, $\overline{h}_a$ & \,\,\,\,\,\,\, $d\cdot t$ &\,\,\,\,\,\,\, $d$ \\ \hline\\ 
\quad moment $n = 2$ \quad \\[2mm]
\quad $t = 16$ h \quad  & $22.25$ & \,\, $3.83$ & $-20.33$ & $-4.14$ & \, $38.69$ & \, $2.4$ &\\
\quad $t = 20$ h \quad  & $24.88$ & \,\, $4.00$ & $-23.03$ & $-4.30$ & \, $35.82$ & \, $1.8$ & \\
\quad $t = 24$ h \quad  & $26.98$ & \,\, $4.11$ & $-24.12$ & $-4.50$ & \, $35.25$ & \, $1.5$ & \\[4mm]
\\ \hline\\
\quad moment $n = 3$ \quad \\[2mm]
\quad $t = 16$ h \quad  & $22.25$ & \,\, $2.59$ & $-20.33$ & $-4.14$ & \, $33.71$  & \, $2.1$ &\\
\quad $t = 20$ h \quad  & $24.88$ & \,\, $2.73$ & $-23.03$ & $-4.30$ & \, $32.31$ & \, $1.6$ & \\
\quad $t = 24$ h \quad  & $26.98$ & \,\, $2.86$ & $-24.12$ & $-4.50$ & \, $29.71$ & \, $1.2$ & \\ \\
\hline\\ 
\quad $Average$  \quad  & $24.70$ & \,\, $3.35$ & $-22.49$ & $-4.31$ & \, $34.25$ & \, $1.8$ & \\ \\
\hline
\end{tabular}
}
\vskip3mm
\caption{Parameter values to fit the second and third moments given in Equation (\ref{eq:moments3}) to those computed from the stochastic simulations. Note that $\overline{h}_a = -\frac{1}{2} \Gamma S(\infty)$, while $d_0 = 0$. \label{table:fitting}}
\end{table}

Equation (\ref{eq:moments}) now gives that 
\begin{eqnarray}
\label{eq:moments1}
\lim_{t \to \infty}\sum_{j=0}^\infty \phi(e^{Kt}q p_0)\frac{m^j e^{-m}}{j!} = { p}^{n(1-\gamma_a-\gamma_b)+\gamma_a(\frac{\beta_a^n-1}{\beta_a-1})+\gamma_b(\frac{\beta_b^n-1}{\beta_b-1})}_\infty,
\end{eqnarray}
assuming that the anastomosis and branching processes are independent,
\[
\beta_b = \overline{h_b}+1=\overline{\alpha}+1=4.35,
\]
and  
\[
\beta_a = \overline{h_a}+1=-\frac{1}{2}\Gamma S(\infty)+1=-3.31,
\]
using the values from Table 1, and since $S(\infty)=a/bc=\sqrt{2K\Gamma+\mu^2}=15.57$.
The upshot of this is that unless we are in the traveling frame of the soliton all the moments are exceedingly small. However, in the 
traveling frame of the soliton $\xi = x-ct+x_0$, the average is just the soliton itself
\[
p= \frac{(2K\Gamma+\mu^2)c}{2\Gamma(c-F_x)}\mbox{sech}^2\left(\frac{\sqrt{2K\Gamma+\mu^2}}{2\Gamma(c-F_x)}
\xi \right),
\]
whereas the second moment is
\[
{p}^{2(1-\gamma_a-\gamma_b)+\gamma_a(\frac{\beta_a^2-1}{\beta_a-1})+\gamma_b(\frac{\beta_b^2-1}{\beta_b-1})}=\left(\frac{(2K\Gamma+\mu^2)c}{2\Gamma(c-F_x)}\right)^{0.59}\mbox{sech}^{1.18}\left(\frac{\sqrt{2K\Gamma+\mu^2}}{2\Gamma(c-F_x)}
\xi \right)
\]
and the nth moment is
\[
p^{\zeta_n}=\left(\frac{(2K\Gamma+\mu^2)c}{2\Gamma(c-F_x)}\right)^{\zeta_n}\mbox{sech}^{2\zeta_n}\left(\frac{\sqrt{2K\Gamma+\mu^2}}{2\Gamma(c-F_x)}
\xi \right), 
\]
where
\[
\zeta_n={n(1-\gamma_a-\gamma_b)+\gamma_a(\frac{\beta_a^n-1}{\beta_a-1})+\gamma_b(\frac{\beta_b^n-1}{\beta_b-1})}.
\]

Now by the argument in the previous section, we get the moments
\begin{eqnarray}
\label{eq:moments3}
\langle p^n \rangle \approx e^{(-nd t-nd_0)}p^{\zeta_n}.
\end{eqnarray}
In particular for the mean, we get a one parameter $(d)$ family,
\[
\langle p \rangle \approx e^{(-d t-d_0)}\frac{(2K\Gamma+\mu^2)c}{2\Gamma(c-F_x)}\mbox{sech}^2\left(\frac{\sqrt{2K\Gamma+\mu^2}}{2\Gamma(c-F_x)}
\xi \right),
\]
whereas for the second moment, we get an expression depending on two parameters, $d$ and $\epsilon$,
\[
\langle p^2 \rangle \approx e^{(-2d t-2d_0)}\left(\frac{(2K\Gamma+\mu^2)c}{2\Gamma(c-F_x)}\right)^{0.6\epsilon}\mbox{sech}^{1.2 \epsilon}\left(\frac{\sqrt{2K\Gamma+\mu^2}}{2\Gamma(c-F_x)}
\xi \right),
\]
where instead of choosing the size of the perturbation in Equation (\ref{eq:angiojumps}) to be $1$, we have chosen it to range between $0$ and $1$, as measured by $\epsilon$, $0 \leq \epsilon \leq 1$.

Comparison of the theory and simulations is shown in Figure 1, for the second and the third moment at increasing times. 
\begin{figure}
\label{fig:sim}
\includegraphics[width=7cm]{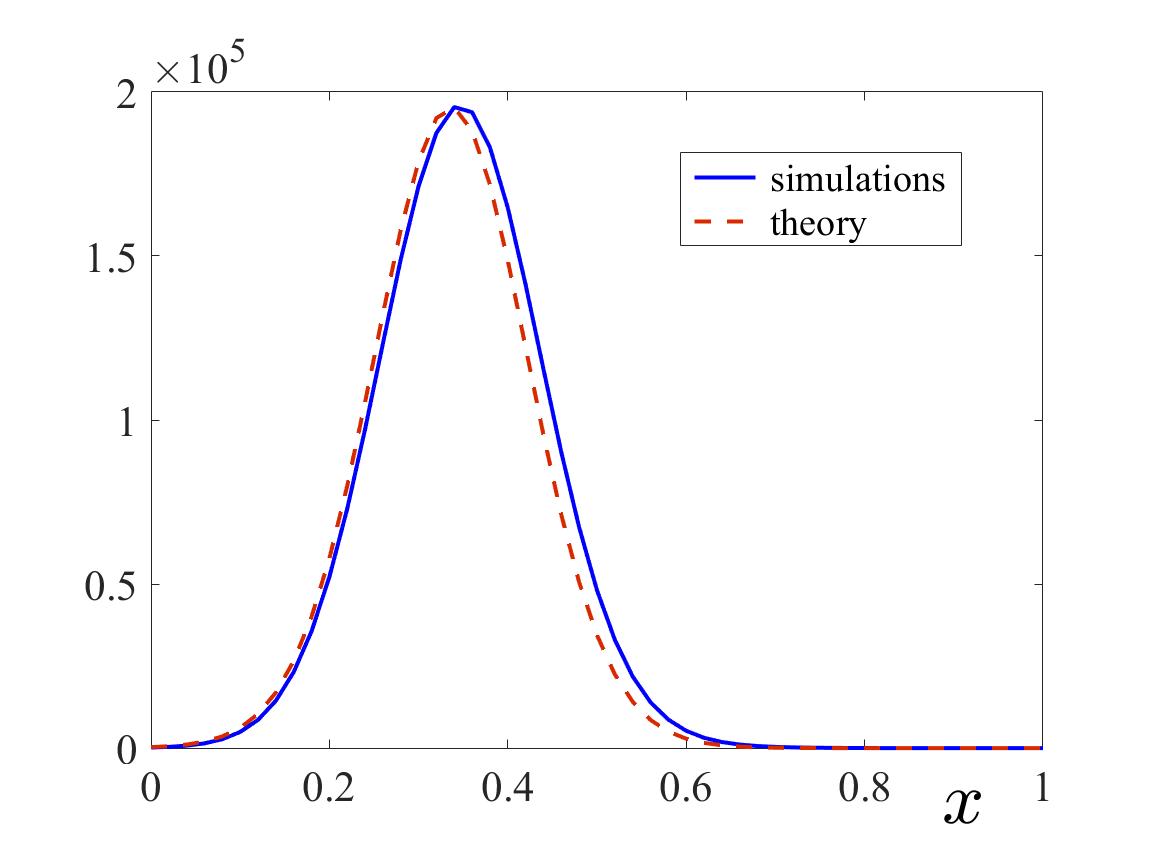}
\includegraphics[width=7cm]{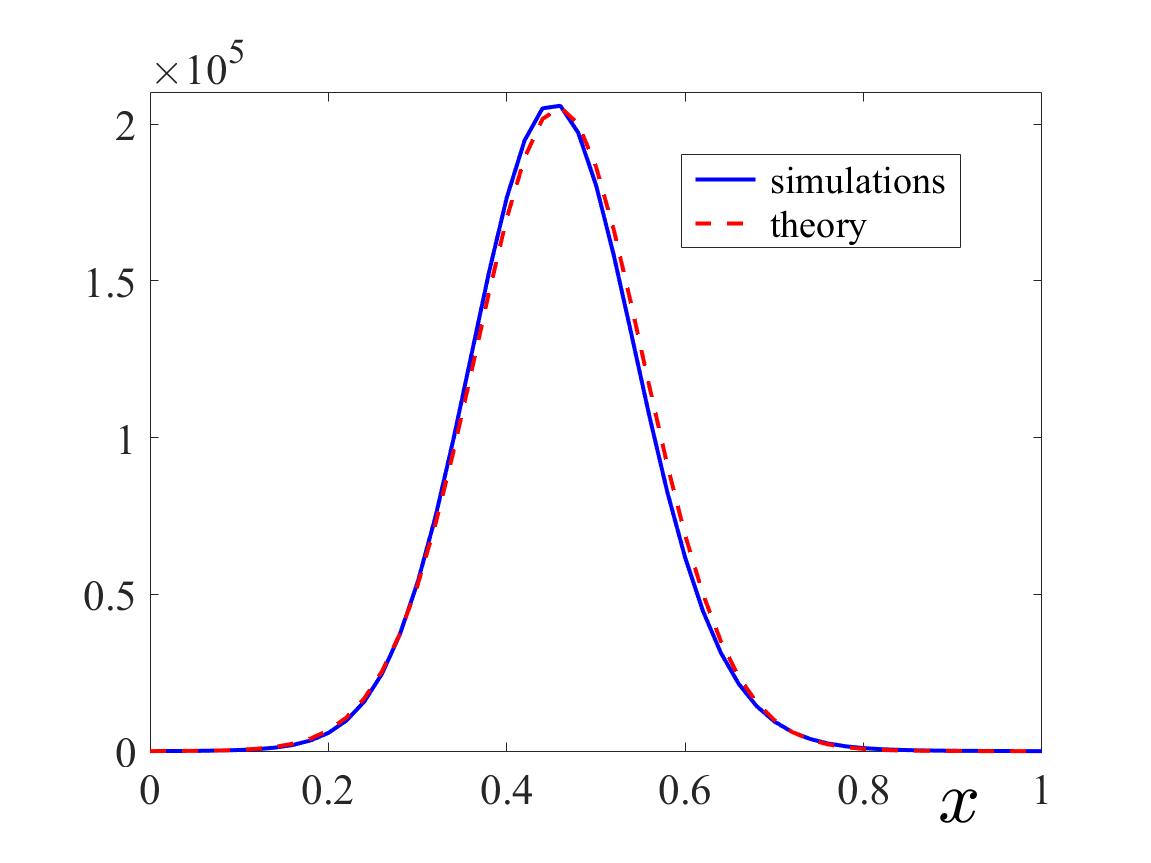}
\includegraphics[width=7cm]{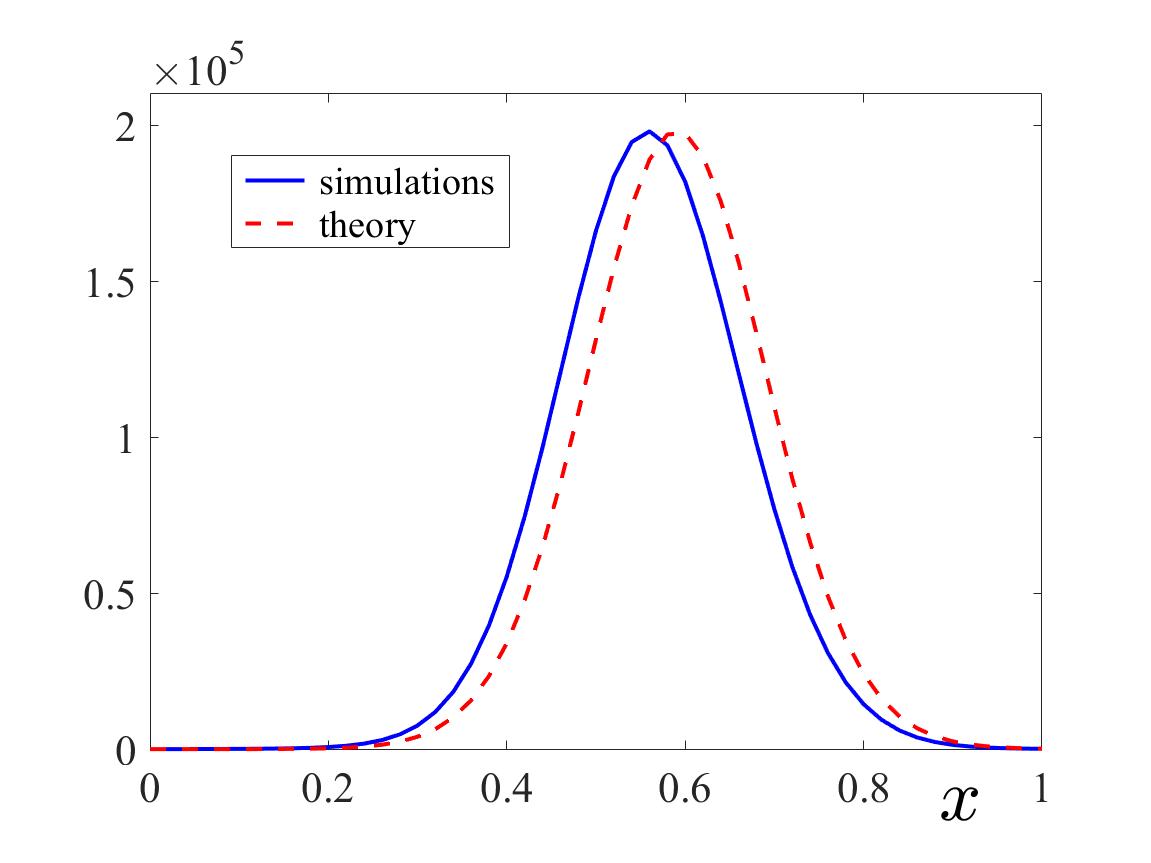}
\includegraphics[width=7cm]{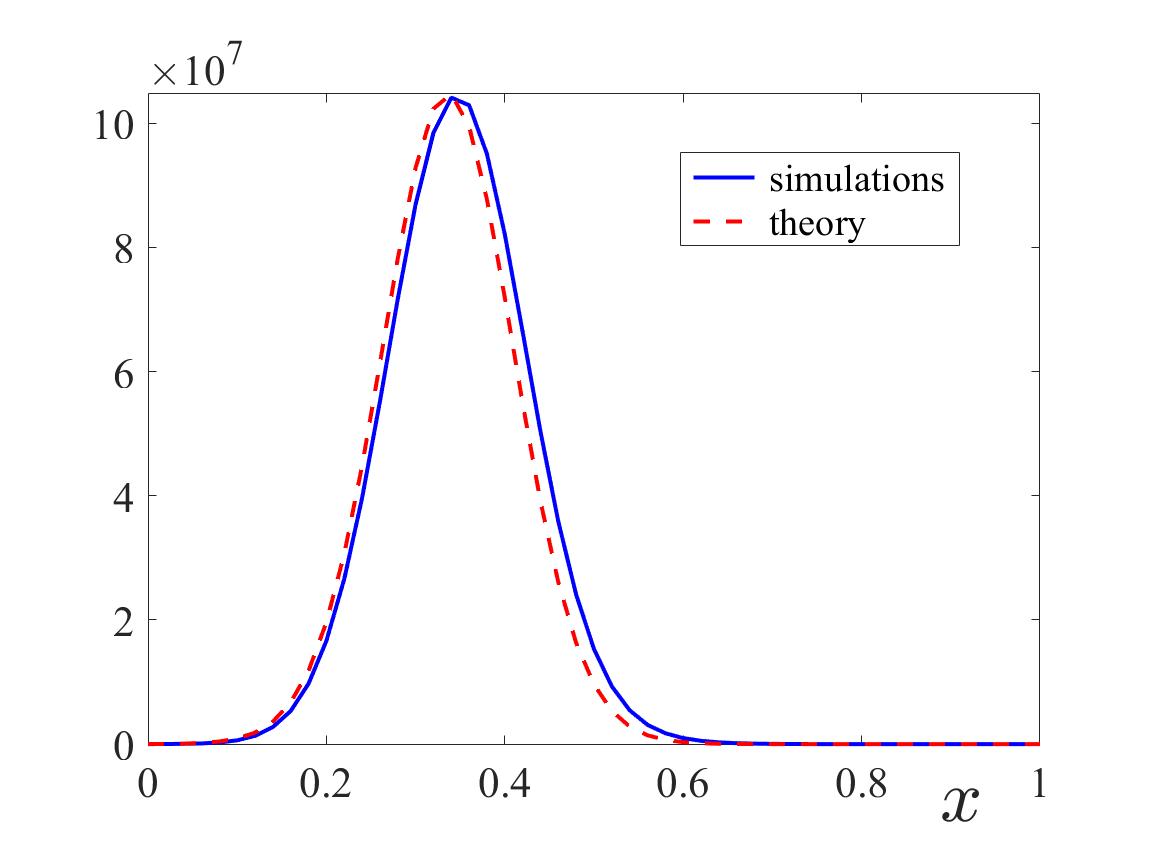}
\includegraphics[width=7cm]{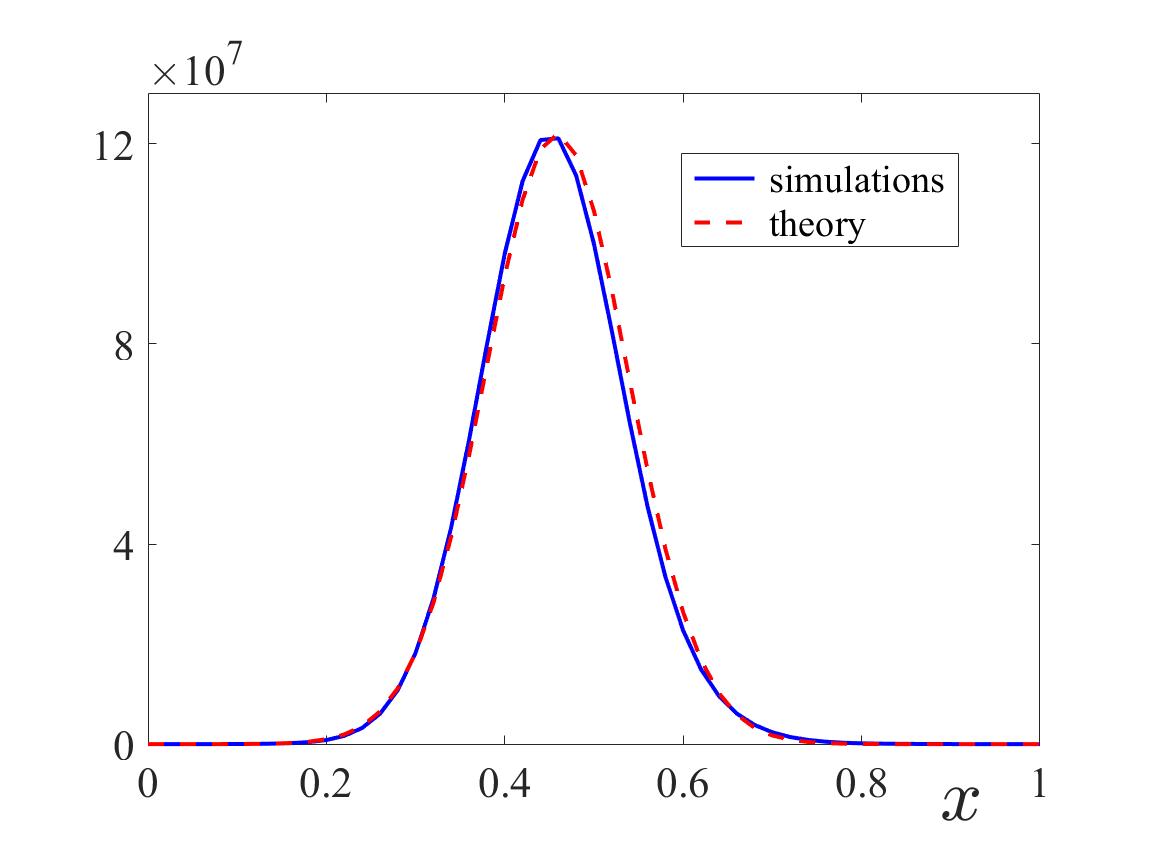}
\includegraphics[width=7cm]{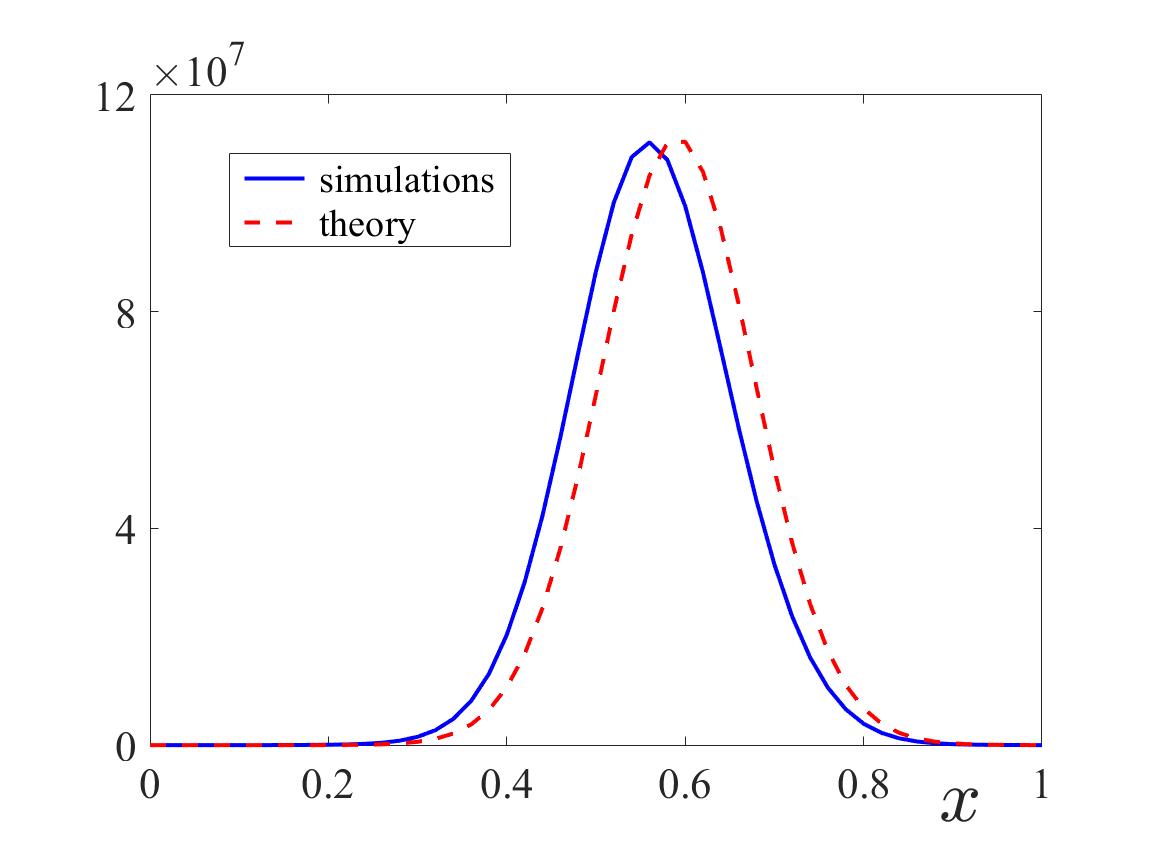}
\caption{A comparison of theory and simulations for the second (first three) and the third (last three) moments of the stochastic angiogenesis equation (\ref{eq:sangiojumps}). The times are 16 hours (left), 20 hours (middle) and 24 hours (right), from left to right.}
\end{figure}

It is reasonable to consider the decaying statistical quantities above instead of taking the limit as $t \to \infty$. We are interested in times until a  blood vessel (artery or vein) encounters another one or joins with a tumor. The decay coefficient is small $d \approx 1.8/hour$ so these vessels shrink slowly as they are elongated. The coefficient $d_0=0$, from the simulations.

There is a small discrepancy between the theoretical moments and the simulated ones at 
24 hours, see Figure 1, both for the second and the third moment. This reason for this is well know. The soliton is not stable but there is a one-dimensional subspace of translated solitons that is stable. This is called orbital stability, see \cite{We87}. The moments are also orbitally stable so the theoretical moments are a small translate of the simulated ones, at 24 hours. This translation distance increases with time.

\section{The Structure Functions}

We compute the moments of the density differences $p_1-p_2=p(\xi_1)-p(\xi_2)$ to probe the fine structure 
in angiogenesis, $\xi_i = x-ct+x_i^0,\ i =1,2$. This means the we are looking at the difference of two solitons that 
differ only in their initial position. 

The density differences $\delta p = p_1-p_2$ is
\[
\delta p =\frac{(2K\Gamma+\mu^2)c}{2\Gamma(c-F_x)}\left (\mbox{sech}^2\left(\frac{\sqrt{2K\Gamma+\mu^2}}{2\Gamma(c-F_x)}
\xi_1\right)-\mbox{sech}^2\left(\frac{\sqrt{2K\Gamma+\mu^2}}{2\Gamma(c-F_x)}
\xi_2\right)\right). 
\]
We let $\xi_1=\xi_0+l$ and $\xi_2=\xi_0$, where $l$ is a lag variable. Then $\delta p$ becomes
\[
\delta p \approx\frac{(2K\Gamma+\mu^2)^{3/2}c}{4\Gamma^2(c-F_x)^2}\left[2\mbox{sech}^2\!\left(\frac{\sqrt{2K\Gamma+\mu^2}}{2\Gamma(c-F_x)}
\xi_0\right) \tanh\!\left(\frac{\sqrt{2K\Gamma+\mu^2}}{2\Gamma(c-F_x)}
\xi_0\right)\right] l. 
\]
by the same arguments as above, the second moment is
\begin{eqnarray*}
&&{\delta p}^{2(1-\gamma_a-\gamma_b)+\gamma_a(\frac{\beta_a^2-1}{\beta_a-1})+\gamma_b(\frac{\beta_b^2-1}{\beta_b-1})}\approx\\
&&\left(\frac{(2K\Gamma+\mu^2)^{3/2}c}{2\Gamma^2(c-F_x)^2}\right)^{0.6}\left[\mbox{sech}^2\!\left(\frac{\sqrt{2K\Gamma+\mu^2}}{2\Gamma(c-F_x)}
\xi_0\right) \tanh\!\left(\frac{\sqrt{2K\Gamma+\mu^2}}{2\Gamma(c-F_x)}
\xi_0\right)\right] ^{0.6} l^{0.6},
\end{eqnarray*}
for $l$ small and the nth moment is
\[
(\delta p)^{\zeta_n}\approx\left(\frac{(2K\Gamma+\mu^2)^{3/2}c}{2\Gamma^2(c-F_x)^2}\right)^{\zeta_n}\left[\mbox{sech}^2\left(\frac{\sqrt{2K\Gamma+\mu^2}}{2\Gamma(c-F_x)}
\xi_0\right)\tanh\left(\frac{\sqrt{2K\Gamma+\mu^2}}{2\Gamma(c-F_x)}
\xi_0\right)\right]^{\zeta_n} l^{\zeta_n}, 
\]
for $l$ small,
where
\[
\zeta_n={n(1-\gamma_a-\gamma_b)+\gamma_a(\frac{\beta_a^n-1}{\beta_a-1})+\gamma_b(\frac{\beta_b^n-1}{\beta_b-1})}.
\]
\vspace{0.5cm}

These formulas then give the structure functions,
\begin{eqnarray}
\label{eq:moments4}
\langle (\delta p)^n \rangle \approx e^{(-nd t-nd_0)}(\delta p)^{\zeta_n},
\end{eqnarray}
as in the previous section.

\section{Discussions}
Considering the scaling exponents of the moments and the structure functions
\[
\zeta_n={n(1-\gamma_a-\gamma_b)+\gamma_a(\frac{\beta_a^n-1}{\beta_a-1})+\gamma_b(\frac{\beta_b^n-1}{\beta_b-1})}, 
\]
we see that the nth moment has a scaling $n(1-\gamma_a-\gamma_b)$, with intermittency corrections $\gamma_a(\frac{\beta_a^n-1}{\beta_a-1})+\gamma_b(\frac{\beta_b^n-1}{\beta_b-1})$. This is different from the three-dimensional stochastic Navier-Stokes equation
\cite{BB211}, where the moments are skewed Gaussians, but the structure functions have a scaling with intermittency corrections. The structure functions of (\ref{eq:sangiojumps}) have the same scaling as the moments, however, moments of (\ref{eq:sangiojumps}) are even functions and reach their maximum at $\xi =0$, whereas the structure functions are odd functions of $\xi_0$ and have two maxima at $\xi_0 =\pm \mbox{sech}^{-1}(\sqrt{3/2})$. 

We thus see that the statistical theory of (\ref{eq:sangiojumps}) is very different from that of the stochastic Navier-Stokes equation. If we are in the traveling frame of the soliton, we see 
decaying soliton-like terms, given by the formulas above. Apart from these, the statistical quantities consist of small and rapidly decaying radiation terms. 

We also see that a much simpler perturbation term, with continuous increments, of the density gives the same results, 
namely
\begin{eqnarray}
\label{eq:sangiojumps_1}
dp =\frac{\sigma}{2}\D_v p+ [\alpha(C) \ \delta_v(v-v_0)- \Gamma  S]\, p- v\cdot \n_x p-\n_v\cdot\!\left[
(F(C)-\sigma v)p\right]dt\no\\
+\eps \left(-(d-\frac{d_0}{2})dt-\sqrt{d_0}db_t \right) p+\eps p\int_{\mathbb{R}}h(z,v,t){\bar N}(dz,dt),
\end{eqnarray}
with $b_0=1$. Thus the (Brownian) noise, in vessel branching, only depends of the density (not all of its Fourier components) and is, in this aspect, similar to the noise in anastomosis, that also only depends on the jumps in the density. 

Finally, the statistical theory of angiogenesis only applies to finite times, until a network of blood vessels is beneficially established in organs and recovering tissue, or malignantly connected a cancerous tumor.

Summarizing, we have proposed a statistical theory of angiogenesis by adding appropriate noise terms to the equation for the density of active tip cells. The shape of the noise terms mimic the birth and death processes of tip branching and anastomosis, respectively, plus some other terms coming from Brownian motion. Thus, our start point is a stochastic PDE with Gaussian and Poisson noises that define an infinite dimensional L\'evy process. We have solved the associated functional equation by finding appropriate integrating factors and therefore obtained the invariant measure. The result is an invariant measure multiplied by an exponentially decaying factor.   By comparing theory and numerical simulations of the underlying angiogenesis model, we have obtained the appropriate values of the parameters involved in the stochastic PDE. We have found that the invariant measure and the moments are functions of the soliton which approximates the deterministic density \cite{BB16,BCT16}.

\bigskip
\noindent
{\bf Acknowledgements}: This work was supported by a Chair of Excellence at the University of Carlos III, in the spring of 2015, that is gratefully acknowledged. LLB, MC and FT's work has been supported by FEDER/Ministerio de Ciencia, Innovaci\'on y Universidades--Agencia Estatal de Investigaci\'on Grant No. PID2020-112796RB-C22, by the Madrid Government (Comunidad de Madrid-Spain) under the Multiannual Agreement with UC3M in the line of Excellence of
University Professors (EPUC3M23), and in the context of the V PRICIT (Regional Programme of Research and Technological Innovation).

\appendix
\setcounter{equation}{0}
\renewcommand{\theequation}{A.\arabic{equation}}
\section{Jumps and L\'evy Processes}
\label{sec:jumps}
In this appendix, we digress to explain the Ito-L\'evy's formula, a generalization of Ito's formula that is necessary to solve (\ref{eq:angiojumps}). 
First we define stochastic processes with jumps, following \O ksendal and Sulem \cite{OS05}, where more information can be found. A L\'evy process is a stochastic process on a filtered
probability space $(\Omega,\cF,\{\cF_t\}_{t\ge 0}, \mathbb{P})$, that takes its values in $\R$ and is continous
in probability and has stationary independent increments. The $\{\cF_t\}_{t\ge 0}$ are an increasing sequence of sigma-algebras indexed by $t$. They are called a filtration on the probability space $(\Omega, \cF, \mathbb{P})$.

Let $b_t$ be a one-dimensional Brownian motion on the probability space $(\Omega, \cF, \mathbb{P})$
and suppose that $N(t,z)$ is the number process of a L\'evy process $\eta_t$. 
An Ito-L\'evy process is a stochastic process on $(\Omega, \cF, \mathbb{P})$ of the form
\begin{eqnarray}
\label{eq:Ito-Levy Process}
dx_t=  u(t,\omega)dt + w(t,\omega)db_t+ \int_\R \gamma(t,z,\omega){\bar N}(dz,dt).
\end{eqnarray}
${\bar N}$ is called the compensated jump measure of $\eta_t$,
\[
{\bar N}(dz,dt)= N(dz,dt)-m(dz)dt, \:\:\: \mathrm{if}\:\:\: |z| < R
\]
and 
\[
{\bar N}(dz,dt)= N(dz,dt), \:\:\: \mathrm{if}\:\:\: |z| \geq R.
\]
where $m(U) =\int_U E(N(dz,1))$ is the so called L\'evy measue of $\eta_t$, see \cite{BB13}, Section 1.5.

The main computational tool in the Ito{'}s calculus for Ito-L\'evy processes is 
the Ito-L\'evy Formula, see \cite{OS05}. It is a generalization of the Ito's formula, so that the jumps are included:
Let $x_t$ be the Ito-L\'evy process,
\[
dx_t=udt+\sqrt{\sigma} db_t+\int_\R \gamma(z,t){\bar N}(dz,dt).
\]
Let $g(x,t)\in C^2([0,\infty) \times \R)$ be twice continuously differentiable. Then
\[
y_t =g(t,x_t)
\]
is also an Ito-L\'evy
process and 
\begin{eqnarray}
\label{eq:Ito-Levy formula}
dy_t&=& \frac{\p g}{\p t}(t,x_t)dt +  \frac{\p g}{\p x}(t,x_t)(udt+\sqrt{\sigma} db_t)+ \frac{\sigma}{2}  \frac{\p^2 g}{\p x^2}(t,x_t)dt\no\\
&+&\int_{z < R}\left(g(t,x_t+\gamma(t,z))-g(t,x_t)-\frac{\p g}{\p x}(t, x_t)\gamma(t,x)\right)m(dz)dt\no\\
&+&\int_\R\left(g(t,x_t+\gamma(t,z))-g(t,x_t)\right){\bar N}(dz,dt).
\end{eqnarray}
Here $(udt+\sqrt{\sigma} db_t)^2$ has been computed by the rules
\[
(dt)^2=dt\cdot db_t = db_t \cdot dt = 0,\:\: \mathrm{but}\:\: (db_t)^2 = dt.
\]

It is illustrative to use the (\ref{eq:Ito-Levy formula}) to solve a differential 
equation. We now do so creating the geometric L\'evy process:
We solve the differential 
equation
\[
dZ_t= Z_t[rdt+\alpha  db_t+\int_\R h(z,t){\bar N}(dz,dt)].
\]
Dividing by $Z_t$
\[
\frac{dZ_t}{Z_t} =rdt+\alpha  db_t+\int_\R h(z,t){\bar N}(dz,dt),
\]
we see that a reasonable guess for the function $g$ is
\[
y_t = \ln(Z_t).
\]
Applying the Ito-L\'evy formula (\ref{eq:Ito-Levy formula}), we get that
\begin{eqnarray*}
dy_t &=& \frac{\p \ln(Z_t)}{\p x}(t,x_t)Z_t(rdt+\alpha  db_t) + \frac{1}{2}  \frac{\p^2 \ln(Z_t)}{\p x^2}(t,x_t)(Z_t)^2(rdt+\alpha  db_t)^2\\
&+&\int_\R\ln(1+h(t,z)){\bar N}(dz,dt)+\int_{\R}\left(\ln(1+h(t,z))-h(t,z)\right)m(dz)dt\\
&=&\frac{Z_t}{Z_t}(rdt+\alpha  db_t) -\frac{1}{2}\frac{(Z_t)^2}{(Z_t)^2}(rdt+\alpha  db_t)^2+\int_\R\ln(1+h(t,z)){\bar N}(dz,dt)\\&+&
\int_{\R}\left(\ln(1+h(t,z))-h(t,z)\right)m(dz)dt\\&=&
(r-\frac{\alpha^2}{2})dt+\alpha db_t+\int_\R\ln(1+h(t,z)){\bar N}(dz,dt)\\&+&
\int_{\R}\left(\ln(1+h(t,z))-h(t,z)\right)m(dz)dt
\end{eqnarray*}
because $\frac{\p g}{\p t} =\frac{\p \ln(z_t)}{\p t}=0$. Integrating 
\begin{eqnarray*}
dy_t =d\ln(Z_t)&=&(r-\frac{\alpha^2}{2})dt+\alpha db_t+
\int_\R\ln(1+h(s,z)){\bar N}(dz,ds)\\&+&
\int_{\R}\left(\ln(1+h(s,z))-h(s,z)\right)m(dz)
\end{eqnarray*}
with respect to $t$ and exponentiating, we get that
\begin{equation}
\label{eq:z}
Z_t=Z_0\ e^{\{(r-\frac{\alpha^2}{2})t+\alpha b_t+\int_0^t
\int_\R\ln(1+h(s,z)){\bar N}(dz,ds)+\int_0^t\int_{\R}\left(\ln(1+h(s,z))-h(s,z)\right)m(dz)ds\}}.
\end{equation}
This process is called the geometric L\'evy process. The first two terms in the exponent correspond to the Ito process and the last two terms (integrals) are contributed by the jump (L\'evy) process.

\setcounter{equation}{0}
\renewcommand{\theequation}{B.\arabic{equation}}
\section{Moments from the stochastic model of angiogenesis}
\label{sec:moments}
In this appendix, we indicate how averages and moments are directly computed by numerical simulation of the stochastic angiogenesis model of Equations  \eqref{eq2} plus random branching and anastomosis processes. As in \cite{ter16}, we run $\mathcal{N}$ realizations of the underlying stochastic process and label each one by $\omega$. Thus, $\omega=1,\ldots,\mathcal{N}$. The `microscopic marginal density' is now
\begin{eqnarray}
\tilde{p}_{\mathcal N}(x,t;\omega)=\sum_{k=1}^{N(t,\omega)}\delta_{\sigma_x}(x-X^k(t,\omega)), \label{eqb1}
\end{eqnarray}
for realization $\omega$, whereas the marginal density $\tilde{p}(t,x)$ is the average of \eqref{eqb1} over all realizations as $\mathcal{N}\to\infty$ and $\sigma_x\to 0$. In practice, we select $\mathcal{N}$ so large that adding more realizations does not change the results of the simulations. Typically $\mathcal{N}=400$ suffices. Numerical values of all other parameters can be found in \cite{ter16}.

To calculate the  $n$ moment of the marginal density, we use a similar formula
\begin{eqnarray}
\label{eq:stochmoments}
\langle \tilde{p}_{\mathcal N}^n \rangle = \frac{1}{\mathcal{N}}\sum_{\omega=1}^\mathcal{N} \left[
\sum_{k=1}^{N(t,\omega)}\delta_{\sigma_x}(x-X^k(t,\omega)) \right]^n,\quad n=1,2,3,\ldots, \label{eqb2}
\end{eqnarray}
in which the ensemble average over realizations is explicitly written. Analogous calculations produce the structure function.

The numerical simulations of this paper consider the simple slab geometry of \cite{ter16,BCT16}. The soliton is calculated using data at a time when its development is completed to solve the collective coordinate equations for $K$ and $c$ in Equation \eqref{eq23}; see \cite{BCT16} for details.

\bibliographystyle{plain}

\end{document}